\documentclass[12pt]{amsart}
\usepackage{amsmath,amscd,amssymb,amsfonts}
\newcommand{\h}{\hbox}
\newcommand{\q}{\quad}
\newcommand{\nin}{\par\noindent}
\newcommand{\bs}{\par\bigskip}
\newcommand{\ms}{\par\medskip}
\newcommand{\sk}{\par\smallskip}
\newcommand{\C}{{\mathbf C}}
\newcommand{\HH}{{\mathbf H}}

\newcommand{\Q}{{\mathbf Q}}
\newcommand{\Z}{{\mathbf Z}}

\newcommand{\OO}{{\mathcal O}}
\newcommand{\Sc}{{\mathcal S}}
\newcommand{\Y}{{\mathcal Y}}
\newcommand{\kk}{\overline{k}}

\newcommand{\pit}{\widetilde{\pi}}
\newcommand{\X}{\widetilde{X}}
\newcommand{\YY}{\widetilde{Y}}
\newcommand{\ZZ}{\widetilde{Z}}
\newcommand{\z}{\widetilde{z}}
\newcommand{\U}{\,\overline{\!U}{}}
\newcommand{\QQ}{\overline{\mathbf Q}}
\newcommand{\p}{{\mathfrak p}}
\newcommand\Hdg{{\rm Hdg}}
\newcommand{\into}{\hookrightarrow}
\newcommand{\simto}{\buildrel\sim\over\rightarrow}
\newcommand{\ssb}{\raise.32pt\h{${\scriptscriptstyle\bullet}$}}
\newcommand{\ssc}{\,\raise.2ex\hbox{${\scriptstyle\circ}$}\,}
\newcommand{\bl}{\bigl}
\newcommand{\br}{\bigr}
\newcommand{\ges}{\geqslant}
\begin{document}
\title{Fields of definition of Hodge loci}
\author{Morihiko Saito}
\address{RIMS Kyoto University, Kyoto 606-8502 Japan}
\email{msaito@kurims.kyoto-u.ac.jp}
\author{Christian Schnell}
\address{Department of Mathematics, Stony Brook University, Stony Brook, NY 11794, USA}
\email{cschnell@math.sunysb.edu}
\begin{abstract}
We show that an irreducible component of the Hodge locus of a polarizable variation of Hodge structure of weight 0 on a smooth complex variety $X$ is defined over an algebraically closed subfield $k$ of finite transcendence degree if $X$ is defined over $k$ and the component contains a $k$-rational point. We also prove a similar assertion for the Hodge locus inside the Hodge bundle if the Hodge bundle together with the connection is defined over $k$. This is closely related with the theory of absolute Hodge classes. The proof uses the spread of the Hodge locus, and is quite similar to the case of the zero locus of an admissible normal function.
\end{abstract}
\maketitle
\centerline{\bf Introduction}
\bs\nin
Let $k$ be an algebraically closed subfield of $\C$ with finite transcendence degree. Let $X$ be a smooth complex algebraic variety defined over $k$. Let $\HH$ be a polarizable variation of $\Z$-Hodge structure of weight 0 on $X$ with underlying $\Z$-local system $\HH_{\Z}$ torsion-free. Let $\Hdg_X(\HH)\subset X$ denote the union of the Hodge loci of all local sections of $\HH_{\Z}$. It is a countable union of closed algebraic subvarieties of $X$ by Cattani, Deligne, and Kaplan \cite{CDK}. In this paper we prove the following.
\ms\nin
{\bf Theorem~1.} {\it Let $Z$ be an irreducible component of $\Hdg_X(\HH)$. If $Z$ contains a $k$-rational point $z_0$ of $X$, then $Z$ is defined over $k$ as a closed subvariety of $X$.}
\ms
The argument is essentially the same as in the case of the zero locus of an admissible normal function \cite{Sa3}, where the spread of closed subvarieties is used in an essential way.
\sk
Let $V$ be the algebraic Hodge bundle over $X$ associated to $\HH$ such that the associated algebraic connection $\nabla:\OO_X(V)\to\Omega_X^1\otimes_{\OO_X}\OO_X(V)$ has regular singularities at infinity, where $\OO_X(V)$ is the sheaf of local sections of $V$ (see \cite{De1}). The last condition is satisfied in the geometric case as is well-known, where the Hodge bundle is given by the relative de Rham cohomology sheaves and the connection by the Gauss-Manin connection. (Note that the algebraic structure on $V$ cannot be obtained correctly without using this regular singularity at infinity.) In this paper we also consider the Hodge locus $\Hdg_V(\HH)$ insider $V$, which is used in the theory of absolute Hodge (de Rham cohomology) classes, see \cite{De4}, \cite{CS}, \cite{Vo2}. It is a countable union of closed algebraic subvarieties of $V$ by \cite{CDK}.
\ms\nin
{\bf Theorem~2.} {\it Assume the Hodge bundle $V$ and the connection $\nabla$ are defined over $k$. Then an irreducible component $Y$ of $\Hdg_V(\HH)$ is defined over $k$ as a closed subvariety of $V$ if it contains a $k$-rational point $y_0$ of $V$.}
\ms
The assumption on the Hodge bundle and the connection is satisfied if these are associated with a smooth projective morphism of complex algebraic varieties $f:\Y\to X$ defined over $k$, where $\HH_{\Z}:=R^{2m}f_*\Z_{\Y}(m)/\h{torsion}$, see \cite{De2} for the Tate twist $(m)$. (Without this Tate twist, the Hodge locus in $V$ cannot be defined over $k$ unless $2\pi i\in k$. This is related with the notion of weakly absolute Hodge class in \cite{Vo2}.)
It seems rather difficult to extend Theorems~1 and 2 to the case where ``an irreducible component" is replaced with ``a connected component" in the statements.
\sk
Let $k_0\subset k$ be a finitely generated subfield such that $X$ is defined over $k_0$.
For $Z$ in Theorem~1, let $k_Z\supset k_0$ be the minimal field of definition of $Z\subset X$ over $k_0$, see \cite[Cor.~4.8.11]{Gro3}. Let $\Sc^{ac}(Z)$ denote the set of algebraically closed subfields $K\subset\C$ such that $K\supset k_0$ and $Z$ has a $K$-rational point. It has a partial ordering by inclusion. We define similarly $k_Y$ and $\Sc^{ac}(Y)$ for $Y$ in Theorem~2. Combining Theorems~1 and 2 with \cite[Cor.~4.8.11]{Gro3}, we get the following.
\sk\nin
{\bf Theorem~3.} {\it In the above notation, there is the minimal element of $\Sc^{ac}(Z)$, and it coincides with the algebraic closure $\kk_Z$ of $k_Z$. A similar assertion holds for $\Sc^{ac}(Y)$ and $k_Y$ if the hypothesis of Theorem~$2$ is satisfied.}
\sk
If we define $\Sc^{fg}(Z)$, $\Sc^{fg}(Y)$ by using finitely generated subfields instead of algebraically closed subfields in the definition of $\Sc^{ac}(Z)$, $\Sc^{ac}(Y)$, then the assertion of Theorem~3 holds with $\kk_Z$, $\kk_Y$, $\Sc^{ac}(Z)$, $\Sc^{ac}(Y)$ replaced respectively by $k_Z$, $k_Y$, $\Sc^{fg}(Z)$, $\Sc^{fg}(Y)$, provided that we consider everything modulo finite extensions (that is, $K$ and $K'$ are identified if $K'$ is a finite extension of $K$).
\sk
Theorem~2 is closely related with the theory of absolute Hodge classes (see \cite{De4}, \cite{Ja}, \cite{CS}, \cite{Vo2}). In this paper we consider only {\it de Rham} cohomology classes as in \cite{CS}, \cite{Vo2} (and not so-called {\it adelic} cohomology classes as in \cite{De4}). Assume $f:\Y\to X$ is defined over $k_0=\Q$ so that $V$ is also defined over $\Q$. (Here $X$ is not necessarily connected, since a $\Q$-variety $X_{\Q}$ satisfying $X=X_{\Q}\times_{{\rm Spec}\,\Q}{\rm Spec}\,\C$ is not necessarily absolutely irreducible.)
Let $V(\QQ)$ denote the set of $\QQ$-rational points of $V$. Set $\Y_{x_0}:=f^{-1}(x_0)$ for a point $x_0\in X$. Then $H^{2m}(\Y_{x_0},\C)$ is canonically identified with $V_{x_0}$.
Let $\alpha$ be a Hodge class in $H^{2m}(\Y_{x_0},\C)=V_{x_0}$, and $Z(\alpha,\QQ)\subset V$ be the closure of $\alpha$ in the $\QQ$-Zariski topology. Let $Y(\alpha)$ be an irreducible component of $\Hdg_V(\HH)$ passing through $\alpha$. From Theorem~2 we can deduce the following.
\ms\nin
{\bf Corollary~1}. {\it In the above notation the following three conditions are equivalent$\,:$
\sk\nin
{\rm (a)} $Y(\alpha)$ is defined over $\QQ$ as a closed subvariety of $V$.
\sk\nin
{\rm (a)$'$} $Y(\alpha)\cap V(\QQ)\ne\emptyset$.
\sk\nin
{\rm (a)$''$} $Z(\alpha,\QQ)\subset Y(\alpha)$.}
\ms
It is conjectured that these equivalent three conditions always hold in the situation of Corollary~1.
(In fact, if the Hodge conjecture is true, then the conditions in Corollary~1 and Corollary~2 below would always hold, see \cite{Vo2} and also Remarks~(2.3) and (3.6)(iv) below.)
We have furthermore the following corollaries which are closely related with results in \cite{Vo2}.
\ms\nin
{\bf Corollary~2}. {\it If the equivalent three conditions in Corollary~$1$ are satisfied, then the following three conditions are equivalent$\,:$
\sk\nin
{\rm (b)} $Y(\alpha)^{\sigma}\subset\Hdg_V(\HH)\,\,(\forall\sigma\in{\rm Gal}(\QQ/\Q))$.
\sk\nin
{\rm (b)$'$} For some $\beta\in Y(\alpha)\cap V(\QQ)$, we have $\beta^{\sigma}\in\Hdg_V(\HH)\,\,(\forall\sigma\in{\rm Gal}(\QQ/\Q))$.
\sk\nin
{\rm (b)$''$} $\alpha$ is absolutely Hodge.}
\ms\nin
{\bf Corollary~3}. {\it Assume $x_0$ is a $\Q$-generic point of $X$. Then we have the equivalences
$${\rm (a)}+{\rm (b)}\Longleftrightarrow{\rm (b)}'\Longleftrightarrow{\rm (b)}''.$$
Moreover, if these equivalent conditions are satisfied, then
\sk\nin
{\rm (c)} $Y(\alpha)$ is finite \'etale over $X$.}
\ms
Note that condition~(b)$'$ logically contains condition~(a)$'$, and $Y(\alpha)^{\sigma}$ for $\sigma\in{\rm Gal}(\QQ/\Q)$ in condition~(b) is meaningful only under condition~(a).
Here $Y(\alpha)^{\sigma}$ denotes the image of $Y(\alpha)$ by the contravariant (or right) action of $\sigma\in {\rm Gal}(\QQ/\Q)$ (which is the inverse of the usual action so that $(Y(\alpha)^{\sigma})^{\sigma'}=Y(\alpha)^{\sigma\sigma'}$).
We can take $\sigma$ from ${\rm Aut}(\C)$ since condition~(a) implies that $Y(\alpha)^{\sigma}$ depends only on the image of $\sigma$ in ${\rm Gal}(\QQ/\Q)$.
Moreover we can take $\sigma$ from $G_{\alpha}\backslash G$ with $G_{\alpha}\subset G:={\rm Gal}(\QQ/\Q)$ the stabilizer of $Y(\alpha)$, and $G_{\alpha}\backslash G$ is finite. Note that the subgroup $G_{\alpha}\subset G$ corresponds to the smallest number field $k_{Y(\alpha)}\subset\QQ$ such that $Y(\alpha)$ is defined over $k_{Y(\alpha)}$ (by using the functoriality of the Galois descent).
In condition~(b)$'$, we can similarly take $\sigma$ from $H_{\beta}\backslash G$ with $H_{\beta}\subset G$ the stabilizer of $\beta$, and $H_{\beta}\backslash G$ is finite.
\sk
The implication ${\rm (a)}+{\rm (b)}\Longrightarrow{\rm (b)}''$ is closely related to Principle~B in \cite{De4} (see also \cite{CS} for the de Rham cohomology class version) in case condition~(c) is satisfied.
In fact, it is also possible to reduce the implication ${\rm (a)}+{\rm (b)}\Longrightarrow{\rm (b)}''$ to Principle~B by using the base change under a certain \'etale morphism defined over $\Q$ in this case (where Theorem~2 is used), see Remark~(3.6)(i) below.
\sk
We thank A.~Tamagawa for answering our question on the reference for the invariance of geometric fundamental groups under base changes by extensions of algebraically closed fields.
We also thank the referee for useful comments.
The first-named author is partially supported by Kakenhi 24540039.
The second-named author is partially supported by NSF grant DMS-1331641, and is very
grateful to Daniel Huybrechts for the opportunity to spend the academic year
2013--2014 at the University of Bonn.  
\sk
In Section~1 we review some basic facts from the theories of variations of Hodge structure and algebraic fundamental groups. In Section~2 we prove the main theorems. In Section~3 we explain the relation with the theory of absolute Hodge classes.
\sk\nin
{\bf Conventions.} In this paper, a variety means a separated scheme of finite type over a field of characteristic 0 although we consider only its closed points. So it is a variety in a classical sense, although it can be reducible or non-reduced in general. However, we need only reduced varieties in this paper (since no non-reduced varieties appear), and the reader may assume that the varieties in this paper are reduced.
\sk
For a complex algebraic variety $X$, we use the classical topology for the constant sheaf $\Z_X$, and the Zariski topology for the structure sheaf $\OO_X$.
\sk
We say that a complex algebraic variety $X$ is defined over a subfield $k$ of $\C$, if there is a variety $X_k$ over $k$ together with an isomorphism $X=X_k\times_{{\rm Spec}\,k}{\rm Spec}\,\C$. If $k$ is a subfield of finite transcendence degree, then $\C$ is sometimes viewed as a universal domain in the sense of \cite{We}.
\bs\bs
\vbox{\centerline{\bf 1. Preliminaries}
\bs\nin
In this section we review some basic facts from the theories of variations of Hodge structure and algebraic fundamental groups.}
\ms\nin
{\bf 1.1.~Variations of Hodge structure on singular varieties.}
Let $Y$ be a complex algebraic variety. Let $L$ be a torsion-free $\Z$-local system on $Y^{\rm an}$, and $F$ be a finite filtration of $\OO_{Y^{\rm an}}\otimes_{\Z}L$ by vector subbundles.
We will say that $\bl(L,(\OO_{Y^{\rm an}}\otimes_{\Z}L,F)\br)$ is a polarizable variation of Hodge structure on $Y$ if so is its pull-back to a desingularization of $Y$.
\sk
Let $\eta$ be a global section of $L$. Assuming $Y$ connected, we have the following
$$\h{$\eta$ has type $(0,0)$ everywhere if it has type $(0,0)$ at one point.}
\leqno(1.1.1)$$
This is reduced to the $Y$ nonsingular case, where the assertion is well-known.
More precisely, it is due to P.~Griffiths \cite[Theorem~7.1]{Gri} in the $Y$ complete case, and
P.~Deligne \cite[Corollary~4.1.2]{De2} in the $\HH$ geometric case. In general, the assertion can be shown by using the existence of a canonical mixed Hodge structure on $H^0(Y,\HH)$ together with the property that the restriction morphism $H^0(Y,\HH)\to\HH_y$ is a morphism of mixed Hodge structure for any $y\in Y$, see \cite{Zu} for the $Y$ curve case (where the last property does not seem to be stated explicitly) and \cite{Sa1} in general. For the proof of (1.1.1), the result in \cite{Zu} may be enough since the assertion can be reduced to the curve case (provided that the reader can verify the above property by himself where he must at least understand the definition of $L^2$ complexes in the Poincar\'e metric).
\ms\nin
{\bf 1.2.~Invariance of geometric fundamental groups by base changes.}
Let $X_k$ be an algebraic variety over an algebraically closed subfield $k$ of $\C$, and $X$ be the base change of $X_k$ by $k\into\C$. There is a canonical isomorphism between the algebraic fundamental groups of $X_k$ and $X$ (see for instance \cite{NTM}).
This assertion follows from the theory of Grothendieck on algebraic fundamental groups by using \cite[Exp.~XIII, Prop.~4.6]{Gro4}, according to A.~Tamagawa.
\sk
By the above invariance together with a well-known comparison theorem for topological and complex algebraic coverings \cite[Exp.~XII, Cor.~5.2]{Gro4}, we get the following:
\sk
Let $\pi:Y\to X$ be a topological finite covering. If $X$ has a structure of complex algebraic variety defined over $k\subset\C$, then $Y$ has a unique structure of complex algebraic variety defined over $k$ which is compatible with the one on $X$ via $\pi$.
\ms\nin
{\bf Remark~1.3.} In the case $X$ is smooth affine, it is rather easy to prove the last assertion except for the uniqueness of the structure over $k$. In fact, a complex algebraic covering $\pi:Y\to X$ can be obtained from an analytic one by applying GAGA to a normal variety which is finite over a smooth compactification of $X$. To get a structure of $k$-variety, consider a finite \'etale morphism $\pi':\Y\to X\times S$ of complex algebraic varieties defined over $k$ such that its restriction over $X\times\{s_0\}$ is isomorphic to $\pi:Y\to X$, where $S$ is a complex affine variety defined over $k$ and $s_0$ is a $k$-generic point of $S$. We may assume that there is a relative smooth compactification of $\Y$ over $S$ such that its boundary is a family of divisors with normal crossings over $S$ (shrinking $S$ if necessary). These are all defined over $k$. Then it is enough to restrict $\pi'$ over $X\times\{s\}$ for some $k$-rational point $s\in S$.
\bs\bs
\vbox{\centerline{\bf 2. Proofs of Theorems 1 and 2}
\bs\nin
In this section we prove the main theorems.}
\ms\nin
{\bf 2.1.~Proof of Theorem~1.} The assertion is reduced to the $X$ affine case by taking a $k$-affine open subset of $X$ containing $z_0$ and by using the closure.
Let $Z^o$ denote the smooth part of $Z$. There is a finite \'etale morphism $\pi_{Z^o}:\ZZ^o\to Z^o$ such that
$$\HH_{\ZZ^o,\Z}:=\pi_{Z^o}^*\HH_{\Z}|_{Z^o}$$
has a global section $\eta_{\ZZ^o}$ which has everywhere type $(0,0)$. Here $\pi_{Z^o}$ is finite by using the polarization restricted to the type $(0,0)$ part (on which the Weil operator is trivial), see \cite{CDK}.
\sk
The morphism $\pi_{Z^o}$ can be extended to a finite morphism
$$\pi_Z:\ZZ\to Z\,\,\,\,\h{with}\,\,\,\,\ZZ\,\,\,\,\h{normal.}$$
This may be obtained, for instance, by considering a smooth variety which is proper over $X$ and contains $\ZZ^o$ as a dense open subvariety, and then using the Stein factorization (i.e., Spec of the direct image of its structure sheaf, see \cite{Ha}).
\sk
Then $\eta_{\ZZ^o}$ is uniquely extended to a global section $\eta_{\ZZ}$ of
$$\HH_{\ZZ,\Z}:=\pi_Z^*\HH_{\Z}|_Z.$$
Indeed, $\eta_{\ZZ^o}$ is a global section defined on the dense Zariski-open subset $\ZZ^o$ of $\ZZ$, and $\ZZ$ is analytic-locally irreducible by using the condition that $\ZZ$ is normal.
\sk
Since $Z$ is affine and $\pi_Z$ is finite, we see that $\ZZ$ is also affine. So there is a closed embedding $\ZZ\into\C^n$, which implies a closed embedding
$$\iota:\ZZ\into \X:=X\times\C^n,$$
such that its composition with the projection $pr_1:\X\to X$ coincides with the composition of $\pi_Z$ with $Z\into X$. Note that $\X$ and the projection $pr_1:\X\to X$ are defined over $k$. Moreover, we may assume that there is a k-rational point $\z_0\in\ZZ$ over $z_0$ by changing the affine coordinates of $\C^n$ if necessary.
(The above argument can be simplified slightly in case the vector bundle $V$ of the variation of Hodge structure $\HH$ is defined over $k$. In fact, $\X$ can be replaced with $V$ which is trivialized by shrinking $X$, and it is unnecessary to assume $\ZZ$ to be normal as long as the section can be extended over $\ZZ$.)
\sk
Consider a spread of $\ZZ\subset\X$ defined over $k$
$$\ZZ_S\subset\X\times S,$$
satisfying
$$\ZZ=\ZZ_{S,s_0}\q\h{in}\,\,\,\,\X,
\leqno(2.1.1)$$
where $S$ is an integral (that is, irreducible and reduced) complex affine variety defined over $k$ with $s_0\in S$ a $k$-generic point, and
$$\ZZ_{S,s}:=\ZZ_S\cap\bl(\X\times\{s\}\br)\subset\X\q\h{for}\,\,\,s\in S,$$
(see for instance \cite{Sa3}). We have
$$z_0\in\ZZ_{S,s}\q\h{for any}\,\,\,s\in S,
\leqno(2.1.2)$$
since
$$\ZZ_S\cap\bl(\{z_0\}\times S\br)=\{z_0\}\times S.$$
Indeed, the left-hand side is defined over $k$ and contains $s_0$.
\sk
There is a sufficiently small open neighborhood $U$ of $s_0$ in classical topology such that
$\eta_{\ZZ}$ can be extended uniquely to a global section $\eta_{\ZZ_U}$ of
$$\HH_{\ZZ_U,\Z}:=\rho_{\ZZ_U}^*\HH_{\Z},$$
where
$$\rho_{\ZZ_U}:\ZZ_U:=\ZZ_S\cap(\X\times U)\to X$$
is the composition of the canonical morphisms $\ZZ_U\to\X\to X$.
This can be shown, for instance, by taking a desingularization of a partial compactification of $\ZZ_S$ which is proper over $S$ and is defined over $k$, where we assume that the total transform of the boundary is a divisor with simple normal crossings over $S$, see also \cite{Sa3}.
(Here it is also possible to apply the generic base change theorem in \cite{De3} together with (1.1.1). One can also use a Whitney stratification defined over $k$, which may have been known to many people for a long time, see \cite{Ma}, \cite{Ve}, \cite{Te}, \cite{Li}, \cite{Ar}, etc., although it seems more difficult to verify all the necessary arguments by oneself. Note that it is much easier to construct controlled vector fields in the normal crossing case than in the general singularity case.)
\sk
The restriction of $\eta_{\ZZ_U}$ to $\ZZ_s$ has everywhere type $(0,0)$ for any $s\in U$ by (1.1.1) and (2.1.2) (In fact, it is a global section of a polarizable variation of Hodge structure, and has type $(0,0)$ at $\z_0$.) This implies
$$Z=pr_1(\ZZ_U),
\leqno(2.1.3)$$
since $Z$ is a union of irreducible components of $\Hdg_X(\HH)$, where the inclusion $\subset$ is clear by (2.1.1). We then get
$$Z=pr_1(\ZZ_S),
\leqno(2.1.4)$$
by considering the pull-backs of local defining functions of $Z$ in $X$ to $\ZZ_S$, which vanish on $\ZZ_U\subset\ZZ_S$.
Since $pr_1(\ZZ_S)$ is defined over $k$, Theorem~1 is proved.
\ms\nin
{\bf 2.2.~Proof of Theorem~2.} Let $Z$ and $z_0$ respectively denote the image of $Y$ and $y_0$ in $X$. By \cite{CDK}, $Y$ is finite over $X$, and $Z$ is a closed subvariety of $X$. We first consider the case
$$\h{$Z$ is an irreducible component of $\Hdg_X(\HH)\subset X$.}
\leqno(2.2.1)$$
(In general we may only have an inclusion of $Z$ in an irreducible component of $\Hdg_X(\HH)$, where it is unclear if $Z$ is defined over $k$.)
By (2.2.1) and Theorem~1, $Z$ is defined over $k$. Take a desingularization $\rho:Z'\to Z$ defined over $k$. Set $V_Z:=V\times_XZ$, and similarly for $V_{Z'}$. We have the commutative diagram
$$\begin{array}{cccccccccccc}
Y&\buildrel{\rho''}\over\longleftarrow&Y'\\
\cap&&\cap\\
V_Z&\buildrel{\rho'}\over\longleftarrow&V_{Z'}\\
\llap{$\vcenter{\hbox{$\scriptstyle p$}}$}\Big\downarrow&&
\Big\downarrow\rlap{$\vcenter{\hbox{$\scriptstyle p'$}}$}\\
Z&\buildrel\rho\over\longleftarrow &Z'\end{array}$$
where $Y'$ is the proper transform of $Y$ so that $Y$ is the image of $Y'$ by $\rho'$.
We have to show that there is a $k$-rational point $y'_0\in Y'$ with $\rho''(y'_0)=y_0$.
(This is non-trivial since only the lower half of the diagram is defined over $k$.)
\sk
Take a general curve $C\subset Z$ defined over $k$ and passing through $z_0:=p(y_0)$. More precisely, this is given by an intersection of general hyperplane sections of $Z$ defined over $k$ and passing through $z_0$. Then $p^{-1}(C)\cap Y$ is a curve passing through $y_0$, since the induced morphism $Y\to Z$ is an isomorphism over its image analytic-locally on $Y$. Let $C'$ be the proper transform of $C$ in $Z'$. We see that the proper transform $Y'$ of $Y$ contains a $k$-rational point given by $(y_0,z'_0)\in V\times_XC'$ if we choose an appropriate $k$-rational point $z'_0\in C'$ over $z_0\in C$.
\sk
By replacing $X$, $V$, and $Y$ respectively with $Z'$, $V\times_XZ'$, and $Y'$, the assertion is then reduced to the case
$$Z=X.
\leqno(2.2.2)$$
Here $Y$ is identified with a {\it multivalued} section of $\HH_{\Z}$ which has everywhere type $(0,0)$.
In particular, it is nonsingular, and is finite \'etale over $X$.
Hence it is defined over $k$ by the theory of Grothendieck on fundamental groups \cite{Gro4}, see (1.2) above.
\sk
We have to show moreover
$$\h{the inclusion $Y\into V$ is defined over $k$.}
\leqno(2.2.3)$$
For this we may assume that $Y$ is a {\it univalued} section of $\HH_{\Z}$ which has everywhere type $(0,0)$ on $X$, by using a commutative diagram as above. (More precisely, using the finite \'etale morphism $\pi:Y\to X$ defined over $k$, we replace $X$ with $Y$, $V$ with its base change by $\pi$, and $Y$ with an irreducible component of its base change by $\pi$.)
\sk
Consider a spread of $Y\subset V$ defined over $k$
$$Y_{S'}\subset V\times S',$$
satisfying
$$Y_{S',s'_0}=Y\q\h{in}\,\,\,\,V,
\leqno(2.2.4)$$
where $S'$ is a complex affine variety defined over $k$ with $s'_0\in S'$ a $k$-generic point, and
$$Y_{S',s'}:=Y_{S'}\cap\bl(V\times\{s'\}\br)\subset V\,\,\,\,\h{for}\,\,\,s'\in S',$$
(see for instance \cite{Sa3}). By the same argument as in the proof of (2.1.2), we have
$$y_o\in Y_{S',s'}\,\,\,\,\h{for any}\,\,\,s'\in S'.
\leqno(2.2.5)$$
\sk
We show that there is a $k$-Zariski open subset
$$U'\subset X\times S'$$
containing $X\times\{s'_0\}$ and such that the restriction of the projection $Y_{S'}\to X\times S'$ over $U'$ is an isomorphism.
Calculating the differential of $Y_{S'}\to X\times S'$, we first get a non-empty $k$-open subvariety $U''\subset Y_{S'}$ which is \'etale over $X\times S'$.
By \cite[Theorem~4.4.3]{Gro2}, we have an open immersion $U''\into\U''$ and a finite morphism $\U''\to X\times S'$, which are defined over $k$ and such that their composition coincides with the canonical morphism $U''\to X\times S'$. Set $D=\U''\setminus U''$. Its image $D'$ in $X\times S'$ is defined over $k$ and hence cannot contain $X\times\{s'_0\}$. (Indeed, if $D'$ contains $X\times\{s'_0\}$, then $D'=X\times S'$ by restricting $D'$ to $\{x\}\times S'$ for any $k$-rational point $x$ of $X$. But this is a contradiction.)
So there is $(x,s'_0)\in X\times\{s'_0\}$ which is not contained in $D'$. This implies that the finite morphism $\U''\to X\times S'$ has degree 1 and
$$D'\cap(X\times\{s'_0\})=\emptyset,$$
(by replacing $U''$ if necessary). So it induces an isomorphism
$$U''\simto U':=(X\times S')\setminus D'.$$
We thus get a section $\eta$ of $\OO_{U'}(V_{U'})$ defined over $k$, where $V_{U'}\to U'$ is the base change of $V\to X$ by the projection $U'\to X$.
\sk
We have the induced relative connection
$$\nabla':\OO_{U'}(V_{U'})\to\Omega^1_{U'/S'}\otimes_{\OO_{U'}}\OO_{U'}(V_{U'}).$$
Consider the zero locus of
$$\nabla'\eta\in\Gamma\bl(U',\Omega^1_{U'/S'}\otimes_{\OO_{U'}}\OO_{U'}(V_{U'})\br).$$
It is a Zariski-closed subset of $U'$, which is defined over $k$ and contains $X\times\{s'_0\}$. Hence it coincides with $U'$. (This is shown by restricting to $\{x\}\times S'$ for any $k$-rational point $x$ of $X$.)
This implies that $Y_{S',s'}$ in (2.2.5) is contained in the image of a horizontal section of $V\to X$ for any $s'\in S'$ (by shrinking $S'$ if necessary). Then (2.2.5) implies that $Y_{S',s'}$ is contained in $Y$. So Theorem~2 is proved in case (2.2.1) is satisfied.
\sk
In case the image of $Y$ in $X$ is a proper subvariety of an irreducible component $Z'$ of $\Hdg_X(\HH)$, we can replace $X$ with a desingularization of $Z'$, and $\HH$ with its pull-back (by using the same argument as above where $C$ is a general curve contained in the image of $Y$ in $X$). Moreover, we can divide $\HH$ by the maximal subvariation of $\Z$-Hodge structure of type $(0,0)$ over $X$ contained in $\HH$. Then we can proceed by induction on the rank of $\HH$. This finishes the proof of Theorem~2.
\ms\nin
{\bf Remark~2.3.} Condition~(a) in Corollary~1 is satisfied if there is a subset $\Sigma\subset Y(\alpha)$ which is not contained in any countable union of proper closed subvarieties of $Y(\alpha)$ and such that the Hodge conjecture holds for any $\beta\in\Sigma$. In fact, for an irreducible component $S$ of the Hilbert scheme of $\Y\to X$ which has a flat family over it whose fibers are closed subschemes of codimension $m$ in fibers of $\Y\to X$, we have the cycle map $S\to V$ over $X$, and the above condition implies that $Y(\alpha)$ coincides with the image of $S\times S'$ for certain two components $S,S'$, where the morphism is given by $cl(\xi_s)-cl(\xi'_{s'})$ for closed points $s,s'$ of $S,S'$ if we denote by $\xi_s$ the associated algebraic cycle over $s$, and similarly for $\xi'_{s'}$. This assertion follows from the countability of the irreducible components of the Hilbert scheme (using the Hilbert polynomials as is well-known, see \cite{Gro1}). Then condition~(a) holds since the Hilbert scheme of $\Y\to X$ is isomorphic to the base change of that for $\Y_k\to X_k$. (The last property follows from the definition of Hilbert schemes using representable functors.) Here we also need the fact that the cycle map can be defined by using relative de Rham cohomology sheaves (or algebraic $D$-modules) so that it is compatible with the base change by $k\into\C$, see \cite{Ja}, \cite{Sa2}.
\sk
It is well-known that condition~(b)$''$ in Corollary~2 is satisfied if the Hodge conjecture holds for $\alpha$, see \cite{De4}. (This can be verified by using a construction as in Remark~(3.6)(ii) below).
\bs\bs
\vbox{\centerline{\bf 3.~Relation with absolute Hodge classes}
\bs\nin
In this section we explain the relation with the theory of absolute Hodge classes.}
\ms\nin
{\bf 3.1.} Corollaries~1, 2 and 3 can be generalized naturally to the case where $\Q$ is replaced by a finitely generated subfield $k$ of $\C$ (so $k$ is not algebraically closed in this section). Here ``absolutely Hodge" means that a Hodge class is ``absolutely Hodge over $k$"; that is, the image of the Hodge class by the action of any element of ${\rm Aut}(\C/k)$ (instead of ${\rm Aut}(\C)={\rm Aut}(\C/\Q)$) is still a Hodge class.
(Here $k$ is countable. In fact, if $k$ is countable, then $k[x]$, $k(x)$ and an algebraic closure of $k$ are countable, where $x$ is an algebraically independent variable.)
\sk
Assume $X$, $\Y$, $V$ in Corollary~1 are base changes of $X_k$, $\Y_k$, $V_k$ by $k\into\C$, and moreover
$$X_k:={\rm Spec}\,A_k,\q V_k:={\rm Spec}\,B_k,$$
by restricting $\Y$, $V$ over an affine open subvariety of $X$ defined over $k$.
Here $A_k$ is a $k$-algebra of finite type, and $B_k={\rm Sym}^{\ssb}_{A_k}M_k^{\vee}$ (the symmetric algebra of the dual $M_k^{\vee}$ of $M_k$ over $A_k$) with
$$M_k:=\Gamma(X_k,R^{2m}f_*\Omega^{\ssb}_{\Y_k/X_k}).$$
We assume that $X_k$ is smooth over $k$. However, it is not necessarily absolutely irreducible, and $X$ is not necessarily connected.
Here $X_k$ is absolutely irreducible over $k':=k(X_k)\cap\kk$ which is a finite extension of $k$.
(In fact, the $k$-variety $X_k$ is actually a variety over $k'\supset k$; that is, the structure morphism $X_k\to {\rm Spec}\,k$ naturally factors through ${\rm Spec}\,k'$, since $\Gamma(X_k,\OO_{X_k})$ contains $k'$ by using the normality of $X_k$. Moreover $k'$ is algebraically closed in the function field of $X_k$. The absolute irreducibility of $X_k$ over $k'$ then follows from the theory of regular extensions by Weil, see, for instance, \cite{La}.
Note that $X_k$ is isomorphic to an irreducible component of $X_{k'}:=X_k\times_{{\rm Spec}\,k}{\rm Spec }\,k'$, although the number of irreducible components of $X_{k'}$ does not necessarily coincide with that of $X$, which is equal to the extension degree $[k':k]$, unless $k'/k$ is Galois.)
So one can assume $X$ connected by replacing $X$ with a connected component if $k$ can be replaced by $k'$.
(This non-connectivity of $X$ is related to the existence of a variety $Z$ over $\C$ which is not homeomorphic to $Z^{\sigma}$ for some $\sigma\in{\rm Aut}(\C)$.)
\sk
Shrinking further $X_k$ if necessary, we have free generators $v_1,\dots,v_r$ of $M_k$ over $A_k$. They induce an isomorphism ($k$-Zariski-locally on $X$)
$$V=X\times\C^r.$$
Assume furthermore
$$A_k=k[x_1,\dots,x_n]/\p,$$
with $\p$ a prime ideal. Then $X$ is identified with the zero locus of $\p$ in $\C^n$ so that each morphism of $k$-algebras $\psi:A_k\to\C$ corresponds to
$$\xi=(\xi_1,\dots,\xi_n)\in X\subset\C^n\,\,\,\h{with}\,\,\,\,\xi_j=\psi(x_j)\,\,(\forall j\in[1,n]).$$
The action of $\sigma\in{\rm Aut}(\C/k)$ on $V$ is expressed by
$$V\ni(\xi_1,\dots,\xi_n;\zeta_1,\dots,\zeta_r)\mapsto(\xi_1^{\sigma},\dots,\xi_n^{\sigma};\zeta_1^{\sigma},\dots,\zeta_r^{\sigma})\in V,$$
where $\xi^{\sigma}$ for $\xi\in\C$ denotes the contravariant (or right) action of $\sigma$ defined by the inverse of the natural action of $\sigma$ on $\C$ so that $(\xi^{\sigma})^{\sigma'}=\xi^{\sigma\sigma'}$ holds.
(If $Z={\rm Spec}\,\C[x]/(f_1,\dots,f_r)$ with $f_i=\sum_{\nu}a_{i,\nu}x^{\nu}$ and $a_{i,\nu}\in\C$, then we have $Z^{\sigma}={\rm Spec}\,\C[x]/(f_1^{\sigma},\dots,f_r^{\sigma})$ with $f_i^{\sigma}=\sum_{\nu}a_{i,\nu}^{\sigma}x^{\nu}$. Here $Z^{\sigma}$ is a scheme which is isomorphic to $Z$ as a scheme, but whose structure over ${\rm Spec}\,\C$ is given by the composition with $\sigma^*:{\rm Spec}\,\C\to{\rm Spec}\,\C$. Hence $Z^{\sigma}$ in this paper means $Z^{\sigma^{-1}}$ in \cite{CS}.)
\sk
Let $\eta\in M_{\C}:=M_k\otimes_k\C$. This can be identified with
$$\h{$(g_1,\dots,g_r)\in A_{\C}^n\,\,\,$ such that $\,\,\,\eta=\sum_i\,g_iv_i$,}$$
where $A_{\C}:=A_k\otimes_k\C$.
The value of the section $\eta$ at $\xi=(\xi_1,\dots,\xi_n)\in X$ is identified with
$$(\xi;g_1(\xi),\dots,g_r(\xi))\in X\times\C^r.$$
Its image by the action of $\sigma\in{\rm Aut}(\C/k)$ is given by
$$(\xi^{\sigma};g_1(\xi)^{\sigma},\dots,g_r(\xi)^{\sigma})\in X\times\C^r.$$
If $\eta\in M_k$, or equivalently, if $g_i\in A_k\,\,(\forall i)$, then the image coincides with
$$(\xi^{\sigma};g_1(\xi^{\sigma}),\dots,g_r(\xi^{\sigma}))\in X\times\C^r,$$
and is still contained in the section $\eta$ (as it should be since $\eta$ is defined over $k$).
This argument together with a construction in Remark~(3.6)(ii) below can be used to show that an algebraic cycle class is absolutely Hodge.
\sk
Let $\alpha\in V$. This defines a prime ideal $\p_{\alpha}\subset B_k$ as in \cite{We}. However, ${\rm Spec}\,B_k/\p_{\alpha}$ is not necessarily absolutely irreducible. Let $Z(\alpha,k)$ denote the zero locus of $\p_{\alpha}$ in $V$. In other words, $Z(\alpha,k)$ is the closure of $\alpha$ in the $k$-Zariski topology of $V$. It is a finite union of irreducible closed subvarieties defined over the algebraic closure $\kk$ of $k$ in $\C$, and we have
\ms\nin
(3.1.1)\q ${\rm Gal}(\kk/k)$ acts transitively on the set of irreducible components of $Z(\alpha,k)$.
\ms\nin
Here the number of irreducible components of $Z(\alpha,k)$ is equal to the extension degree $[k_{\alpha}:k]$ with $k_{\alpha}$ the algebraic closure of $k$ in the field of fractions of $B_k/\p_{\alpha}$ (by using the theory of regular extensions as is explained in the beginning of this section).
\sk
Let $Z(\alpha,\kk)\subset V$ denote the zero locus of the ideal of $B_k\otimes_k\kk$ associated with $\alpha$. This is the smallest irreducible closed subvariety of $V$ containing $\alpha$ and defined over $\kk$. In other words, $Z(\alpha,\kk)$ is the closure of $\alpha$ in the $\kk$-Zariski topology. We have moreover
$$\h{$Z(\alpha,\kk)$ coincides with the irreducible component of $Z(\alpha,k)$ containing $\alpha$.}
\leqno(3.1.2)$$
Here the uniqueness of the component follows from a well-known assertion that the singular locus of $Z(\alpha,k)$ is defined over $k$.
Note that, if $G'_{\alpha}\subset{\rm Gal}(\kk/k)$ denotes the stabilizer of $Z(\alpha,\kk)$, then the subgroup $G'_{\alpha}$ corresponds to the finite extension $k_{\alpha}\supset k$, where $k_{\alpha}$ is as in a remark after (3.1.1). Moreover this corresponding field $k_{\alpha}$ coincides with the minimal field of definition of $Z(\alpha,\kk)$ in $V$ by using the functoriality of the Galois descent (which we apply to a finite Galois extension of $k_{\alpha}$ over which $Z(\alpha,\kk)$ is defined).
\sk
Let $x_0$ be the image of $\alpha$ in $X$. If this is a $k$-generic point of $X$, then it defines the zero ideal of $A_k$, and ${\rm Spec}\,B_k/\p_{\alpha}$ is dominant over ${\rm Spec}\,A_k$. So we get
$$\dim Z(\alpha,k)\ges\dim X\,\,\,\h{if $x_0$ is a $k$-generic point.}
\leqno(3.1.3)$$
\sk
We have the following well-known assertion, which shows that it is actually enough to consider all the embeddings $K\into\C$ (instead of all the automorphisms of $\C$) over $k$ in order to see whether $\alpha$ is absolutely Hodge over $k$, where $K\subset\C$ is a fixed subfield containing $k$ and over which $\alpha$ is defined, see also \cite{Vo2}, \cite{CS}.
\ms\nin
{\bf Proposition~3.2.} {\it The set of points of $V$ which are conjugate to $\alpha$ by automorphisms of $\C$ over $k$ is the complement in $Z(\alpha,k)$ of the countable union of all the proper closed subvarieties defined over $k$.}
\ms\nin
{\it Proof.} Each point $\beta$ of $ Z(\alpha,k)$ corresponds to a morphism of $k$-algebras
$$\psi_{\beta}:B_k/\p_{\alpha}\to\C.$$
This is not necessarily injective. It holds if and only if $\beta$ is not contained in any proper closed subvariety of $Z(\alpha,k)$ defined over $k$.
\sk
By the definition of $\p_{\alpha}$, we have the injectivity of
$$\psi_{\alpha}:B_k/\p_{\alpha}\into\C,$$
and $\alpha^{\sigma}\in V$ for $\sigma\in{\rm Aut}(\C/k)$ corresponds to the composition
$$B_k/\p_{\alpha}\buildrel{\psi_{\alpha}}\over\into\C\buildrel{\sigma^{-1}}\over\longrightarrow\C,$$
where $\sigma^{-1}$ appears in order to get the contravariant action. Note that $\p_{\alpha}=\p_{\alpha^{\sigma}}$ since $B_k$ is a quotient of a polynomial ring over $k$.
\sk
Let $K$ be the field of fractions of $B_k/\p_{\alpha}$. Proposition~(3.2) is then reduced to the following well-known assertion:
\sk\nin
(3.2.1)\q For any two inclusions $\psi_a:K\into\C$ over $k$ ($a=1,2)$, there is $\sigma\in{\rm Aut}(\C/k)$ with $\sigma\ssc\psi_1=\psi_2$.
\sk\nin
(This can be shown by using the theory of transcendence bases.)
This finishes the proof of Proposition~(3.2).
\sk
We prove Corollaries~1, 2 and 3 with $\Q$ replaced by $k$, where ``absolutely Hodge" means ``absolutely Hodge over $k$". The reader may assume $k=\Q$.
\ms\nin
{\bf 3.3.~Proof of Corollary~1.}
The implications (a) $\Longrightarrow$ (a)$''$ $\Longrightarrow$ (a)$'$ are clear, and we have (a)$'$ $\Longrightarrow$ (a) by Theorem~2. So Corollary~1 is proved.
\ms\nin
{\bf 3.4.~Proof of Corollary~2.}
Conditions~(b) clearly implies (b)$''$ by using condition~(a)$''$.
Assume condition~(b)$''$ holds, i.e., $\alpha$ is absolutely Hodge over $k$. By Proposition~(3.2) we have
$$Z(\alpha,k)\subset\Hdg_V(\HH),
\leqno(3.4.1)$$
by using the fact that $\Hdg_V(\HH)$ is locally a finite union of closed analytic subsets (since it is contained in the $\Z$-lattice). Note that (3.4.1) follows from (b)$''$ without assuming (a).
Combining (3.4.1) with condition~(a)$''$, we get condition~(b)$'$.
\sk
Assume condition~(b)$'$. We see that $Y(\alpha)^{\sigma}$ in condition~(a)$''$ is contained in the vector subbundle corresponding to the Hodge filtration $F^0\OO_X(V)$ since the latter is defined over $k$. We have to show that $Y(\alpha)^{\sigma}$ is also contained in the analytic subset of $V$ corresponding to $\HH_{\Z}$. Here $Y(\alpha)^{\sigma}$ can be identified with a flat section of $V$ locally on $X$ in $k$-\'etale topology by using the fact that the Gauss-Manin connection is defined over $k$. In fact, let $\YY(\alpha)$ be the normalization of the $k$-closure of $Y(\alpha)$ in $V$, which is identified with the disjoint union
$$\h{$\bigsqcup_{\,i=1}^{\,d}\,Y(\alpha)^{\sigma_i}$},$$
where the $\sigma_i$ are representatives of $G_{\alpha}\backslash G$ with $G_{\alpha}\subset G:={\rm Gal}(\kk/k)$ the stabilizer of $Y(\alpha)$, and $d:=|G_{\alpha}\backslash G|$. There is a canonical \'etale morphism defined over $k$
$$\pit:\YY(\alpha)\to X,$$
which factors through
$$\h{$\bigcup_iY(\alpha)^{\sigma_i}\subset V$}.$$
This defines a section $\eta$ of the vector bundle
$$V\times_X\YY(\alpha)\to\YY(\alpha),$$
which is the base change of $V\to X$ by $\pit$. By definition the image of the section $\eta$ in $V$ coincides with $\bigcup_iY(\alpha)^{\sigma_i}$. This section $\eta$ is flat, since it is defined over $k$ and is flat over the connected component
$$Y(\alpha)=Y(\alpha)^{\sigma_1}\subset\YY(\alpha),$$
where $\sigma_1=id$. So we get the desired flatness of $Y(\alpha)^{\sigma}$ since $Y(\alpha)^{\sigma}=Y(\alpha)^{\sigma_i}$ for some $i$.
This implies that $Y(\alpha)^{\sigma}$ is contained in the analytic subset of $V$ corresponding to $\HH_{\Z}$, since $Y(\alpha)^{\sigma}$ is contained in the subset at $\beta^{\sigma}\in Y(\alpha)^{\sigma}$ by condition~(b)$'$.
So condition~(b) holds. Thus Corollary~2 is proved.
\ms\nin
{\bf 3.5.~Proof of Corollary~3.}
By Corollaries~1 and 2, it remains to show (b)$''\Longrightarrow{\rm (a)}''+{\rm (c)}$.
Assume condition~(b)$''$. Then (3.4.1) holds. Together with (3.1.3) this implies
$$\dim Z(\alpha,k)=\dim X=\dim Y(\alpha).$$
Hence $Y(\alpha)$ is finite \'etale over $X$, and condition~(c) holds. This implies that $Y(\alpha)$ is the unique irreducible component of $\Hdg_V(\HH)$ containing $\alpha$, and we get
$$Z(\alpha,\kk)=Y(\alpha).
\leqno(3.5.1)$$
So condition~(a)$''$ also holds. Thus Corollary~3 is proved.
\ms\nin
{\bf Remarks~3.6.} (i) We can also prove the implication (b)$'$ $\Longrightarrow$ (b)$''$ by using Principle~B (see Remark~(iii) below) where the base change by the \'etale morphism $\pit:\YY(\alpha)\to X$ as in (3.3) is also used. Here we need Theorem~2 to show that $Y(\alpha)$ is defined over $\kk$.
\sk
(ii) If $\alpha\in V$ is associated with an algebraic cycle and if the cycle is extended to a family of algebraic cycles over $X$ defined over $k$, then it corresponds to an element of $M_k$ which corresponds to a section of $V\to X$ defined over $k$. The second condition is satisfied if we take a finitely generated $k$-subalgebra $A'_k$ of $\C$ over which $\Y_{x_0}$ and the cycles are defined, and if we replace $f:\Y\to X$ by the family over ${\rm Spec}\,A'_k\otimes_k\C$, see \cite{Vo1}, \cite{Sa2}. (Here ${\rm Spec}\,A'_k$ is not necessarily absolutely irreducible, and ${\rm Spec}\,A'_k\otimes_k\C$ is not necessarily connected.) The latter construction can be used to show that an algebraic cycle class is absolutely Hodge.
\sk
(iii) Let $f:\Y\to X$ and $V$ be as in Corollary~1. Let $\eta$ be a flat section of $V$ defined on a connected component of $X$. Assume $X_k$ is irreducible. Then Principle~B in \cite{De4}, \cite{CS} asserts that, if the value of $\eta$ at one point is absolutely Hodge over $k$, then the value at any point is. This follows from the assertion (1.1.1) which is a consequence of the global invariant cycle theorem \cite{De2} in this case. Indeed, (1.1.1) together with the assumption of Principle~B implies that the image $\eta^{\sigma}$ of the section $\eta$ by the action of $\sigma\in{\rm Aut}(\C/k)$ has type $(0,0)$ at every point.
\sk
(iv) If there is a subset $\Sigma\subset Y(\alpha)$ which is not contained in any countable union of proper closed subvarieties of $Y(\alpha)$ and such that any $\beta\in\Sigma$ is absolutely Hodge, then conditions~(a) and (b) would hold in the notation of Corollaries~1 and 2, where $\Q$ is replaced with $k$ as in (3.1). In fact, consider all the closed subvarieties of $V$ defined over $\kk$ and not entirely containing $Y(\alpha)$. These are countable, and there is $\beta\in \Sigma$ which is not contained in any such varieties, and is not a singular point of $\Hdg_V(\HH)$. Then the $\kk$-Zariski closure of $\beta$ in $V$ is contained in $Y(\alpha)$ and its images by the action of ${\rm Gal}(\kk/k)$ are contained in $\Hdg_V(\HH)$, since $\beta$ is absolutely Hodge over $k$. (Here the condition that $\beta$ is a smooth point of $\Hdg_V(\HH)$ implies that $Y(\alpha)$ is the only irreducible component containing it.) If this $\kk$-closure does not coincide with $Y(\alpha)$, then this contradicts the condition that $\beta$ is not contained in any closed subvariety of $V$ which is defined over $\kk$ and does not contain entirely $Y(\alpha)$. So conditions~(a) and (b) follow. (See also \cite{Vo2}.)


\begin{thebibliography}{NaTaMo}
\bibitem[Ar]{Ar} Arapura, D., An abelian category of motivic sheaves, Adv.\ Math.\ 233 (2013), 135--195.
\bibitem[CaDeKa]{CDK} Cattani, E., Deligne, P.\ and Kaplan, A., On the locus of Hodge classes, J.\ Amer.\ Math.\ Soc.\ 8 (1995), 483--506.
\bibitem[ChSch]{CS} Charles, F.\ and Schnell, Ch., Notes on absolute Hodge classes (preprint).
\bibitem[De1]{De1} Deligne, P., Equations diff\'erentielles \`a points singuliers r\'eguliers, Lect.\ Notes in Math.\ 163, Springer, Berlin, 1970.
\bibitem[De2]{De2} Deligne, P., Th\'eorie de Hodge II, Publ.\ Math.\ IHES 40 (1971), 5--58.
\bibitem[De3]{De3} Deligne, P., Th\'eor\`eme de finitude en cohomologie $\ell$-adique, Lect.\ Notes in Math.\ 569, Springer, Berlin, 1977, pp.~233--261.
\bibitem[De4]{De4} Deligne, P., Hodge cycles on abelian varieties (notes by J.\ S.\ Milne), Lect.\ Notes in Math.\ 900, Springer, Berlin, 1982, pp.~9--100.
\bibitem[Gri]{Gri} Griffiths, P.~A., Periods of integrals on algebraic manifolds III. Some global differential-geometric properties of the period mapping, Publ.\ Math.\ IHES 38 (1970), 125--180.
\bibitem[Gro1]{Gro1} Grothendieck, A., Techniques de construction et th\'eor\`emes d'existence en g\'eom\'etrie alg\'ebrique IV: les sch\'emas de Hilbert, S\'eminaire Bourbaki, 1960--1961, Exp.\ 221, pp.~249--276.
\bibitem[Gro2]{Gro2} Grothendieck, A., El\'ements de g\'eom\'etrie alg\'ebrique, III-1, Publ.\ Math.\ IHES 11, 1961.
\bibitem[Gro3]{Gro3} Grothendieck, A., El\'ements de g\'eom\'etrie alg\'ebrique, IV-2, Publ.\ Math.\ IHES 24, 1965.
\bibitem[Gro4]{Gro4} Grothendieck, A., Rev\^etements \'etales et groupe fondamental, S\'eminaire de G\'eom\'etrie Alg\'ebrique 1, Lect.\ notes in Math.\ 224, Springer, Berlin, 1971.
\bibitem[Ha]{Ha} Hartshorne, R., Algebraic Geometry, Springer, Berlin, 1977.
\bibitem[Ja]{Ja} Jannsen, U., Mixed motives and algebraic K-theory, Lect.\ Notes in Math.\ 1400, Springer, Berlin, 1990.
\bibitem[La]{La} Lang, S., Algebra, Addison-Wesley, 1993.
\bibitem[Li]{Li} Lipman, J., Equisingularity and simultaneous resolution of singularities, in Resolution of singularities (Obergurgl, 1997), Progr.\ Math.\ 181, Birkh\"auser, Basel, 2000, pp.~485--505.
\bibitem[Ma]{Ma} Mather, J., Notes on topological stability (Harvard, 1970), Bull.\ Amer.\ Math.\ Soc.\ (N.S.) 49 (2012), 475--506.
\bibitem[NaTaMo]{NTM} Nakamura, H., Tamagawa, A.\ and Mochizuki, S., The Grothendieck conjecture on the fundamental groups of algebraic curves, Sugaku Expositions 14 (2001), 31--53.
\bibitem[Sa1]{Sa1} Saito, M., Mixed Hodge modules, Publ.\ RIMS, Kyoto Univ.\ 26 (1990), 221--333.
\bibitem[Sa2]{Sa2} Saito, M., Arithmetic mixed sheaves, Inv.\ Math.\ 144 (2001), 533--569.
\bibitem[Sa3]{Sa3} Saito, M., Normal functions and spread of zero locus (arXiv:1304.3923).
\bibitem[Te]{Te} Teissier, B., Vari\'et\'es polaires, II, Multiplicit\'es polaires, sections planes, et conditions de Whitney, Lect.\ Notes in Math.\ 961, Springer, Berlin, 1982, pp.~314--491.
\bibitem[Ve]{Ve} Verdier, J.-L., Stratifications de Whitney et th\'eor\`eme de Bertini-Sard, Inv.\ Math.\ 36 (1976), 295--312.
\bibitem[Vo1]{Vo1} Voisin, C., Transcendental methods in the study of algebraic cycles, Lect.\ Notes in Math.\ 1594, Sringer, Berlin, 1994, pp. 153--222.
\bibitem[Vo2]{Vo2} Voisin, C., Hodge loci and absolute Hodge classes, Compos.\ Math.\ 143 (2007), 945--958.
\bibitem[We]{We} Weil, A., Foundation of algebraic geometry, AMS 1946.
\bibitem[Zu]{Zu} Zucker, S., Hodge theory with degenerating coefficients. $L_2$ cohomology in the Poincar\'e metric, Ann.\ Math.\ (2) 109 (1979), 415--476.
\end{thebibliography}
\end{document}